\documentclass{article}
\usepackage{amsmath,amssymb,theorem}
\usepackage[french]{babel}

\newcommand{\nobracket}{}
\newcommand{\nocomma}{}
\newcommand{\noplus}{}

\newcommand{\tmop}[1]{\ensuremath{\operatorname{#1}}}
\newcommand{\tmsamp}[1]{\textsf{#1}}
\newcommand{\tmstrong}[1]{\textbf{#1}}
\newcommand{\tmtextbf}[1]{{\bfseries{#1}}}
\newenvironment{proof}{\noindent\textbf{Preuve\ }}{\hspace*{\fill}$\Box$\medskip}
\newtheorem{lemma}{Lemme}
\newtheorem{proposition}{Proposition}
{\theorembodyfont{\rmfamily}\newtheorem{remark}{Remarque}}
\newtheorem{theorem}{Th\'eor\`eme}
\begin{document}

\title{R{\'e}solution du $\partial \bar{\partial}$ pour les courants
prolongeables d{\'e}finis sur un domaine fortement pseudoconvexe d'une
vari{\'e}t{\'e} contractile}
\author{Eramane Bodian, Dian Diallo et Salomon Sambou}
\maketitle

\begin{abstract}
  On r{\'e}sout le $\partial \bar{\partial}$ pour les courants prolongeables
  d{\'e}finis dans un domaine strictement pseudoconvexe d'une vari{\'e}t{\'e}
  contractile.
  
  \ \ \ \ \ \ \ \ \ \ \ \ \ \ \ \ \ \ \ \ \ \ \ \ \ \ \ \ \ \ \ \ \ \ \ \ \ \
  \ {\textbf{Abstract}}
  
  We solve the $\partial \bar{\partial}$-problem for extensible currents
  defined on a strongly pseudoconvex domain of a contractible manifold.
  
  {\tmstrong{Mots clefs}} : Courant prolongeable , $\partial \bar{\partial}$ ,
  cohomologie de De Rham.
  
  {\tmsamp{Classification math{\'e}matique }}2010 : 32F32.
\end{abstract}

\section{Introduction}

Soit $\Omega \Subset X$ un ouvert , on se pose la question suivante :

Si $C$ est un courant $d$-ferm{\'e} sur $\Omega$ , existe - t - il un courant
prolongeable sur $\Omega$ tel que $\partial \bar{\partial} u = C$ ?

Tenant compte de consid{\'e}rations classiques , nous devons pour r{\'e}pondre
{\`a} cette question , avoir {\`a} r{\'e}soudre une {\'e}quation $(\ast) du =
C$ , o{\`u} $C$ est un courant prolongeable, la solution obtenue se
d{\'e}compose en une partie $\partial$-ferm{\'e}e et l'autre $\bar{\partial}$-
ferm{\'e}e. Il est n{\'e}cessaire d'avoir des conditions g{\'e}om{\'e}triques
sur $\Omega$ pour obtenir des solutions du $\partial$ respectivement du
$\bar{\partial}$ pour les courants prolongeables. Partant de r{\'e}sultats
connus de cohomologie de De Rham , nous montrons que si

$\Omega$ est un domaine strictement pseudoconvexe, alors l'{\'e}quation
$(\ast)$ admet une solution sous certaines hypoth{\`e}ses cohomologiques. La
r{\'e}solution du $\partial \bar{\partial}$ devient alors une cons{\'e}quence
des r{\'e}sultats de r{\'e}solution du $\bar{\partial}$ pour les courants
prolongeables obtenus dans \cite{Samb}. Pour finir, on donnera
quelques exemples de domaines $\Omega$.

\section{Pr{\'e}liminaires et notations}

\subsection{D{\'e}finition}

Soit $X$ une vari{\'e}t{\'e} diff{\'e}rentiable , $\Omega \subset X$ un
domaine. Un courant $C$ d{\'e}fini sur $\Omega$ est dit prolongeable s'il
existe un courant $\tilde{C}$ d{\'e}fini sur $X$ tel que $\tilde{C}_{| \Omega}
= C$.

Les op{\'e}rateurs $\partial$ et $\bar{\partial}$ d{\'e}finis pour les formes
diff{\'e}rentielles s'{\'e}tendent aux courants par dualit{\'e}.
\[ \langle \bar{\partial} C, \varphi \rangle = (- 1)^{\deg C + 1} \langle C,
   \bar{\varphi} \rangle . \]
D'apr{\`e}s Martineau \cite{Mart} , si \ $\dot{\bar{\Omega}} =
\Omega$ alors les courants prolongeables de degr{\'e} $p$ sur $\Omega$ sont
{\'e}gaux au dual topologique des $(n - p)$-formes diff{\'e}rentielles de
classe $\mathcal{C}^{\infty}$ sur $X$ {\`a} support sur $\bar{\Omega}$ ,
o{\`u} $n$ d{\'e}signe la dimension de $X$. On note $\check{D}_X^p
(\Omega)$les $p$-courants d{\'e}finis sur $\Omega$ et prolongeables {\`a} $X$,
$D^p ( \bar{\Omega})$ les $p$-formes diff{\'e}rentielles de classe
$\mathcal{C}^{\infty}$ sur $X$ {\`a} support dans $\bar{\Omega}$. Si $X$ est
une vari{\'e}t{\'e} complexe de dimension $n$ , on note $\check{D}_X^{p, q}
(\Omega)$ l'espace des $(p, q)$-courants prolongeables d{\'e}finis sur
$\Omega$ et $D^{p, q} ( \bar{\Omega})$ l'espace des $(p, q)$-formes
diff{\'e}rentielles {\`a} support $\bar{\Omega}$. On note $\check{H}^p
(\Omega)$ le $p^{\tmop{iem}}$ groupe de cohomologie de De Rham des courants
prolongeables d{\'e}finis sur $\Omega$ , $\check{H}^{p, q} (\Omega)$ le $(p,
q)^{\tmop{iem}}$ groupe de cohomologie de Dolbeault des courants prolongeables
d{\'e}finis sur $\Omega$. Si $F \subset X$, alors $H_{\infty}^p (F)$
d{\'e}signe le $p^{\tmop{iem}}$ groupe de cohomologie de De Rham des
$p$-formes diff{\'e}rentielles de classe $\mathcal{C}^{\infty}$ d{\'e}finis
sur $X$, $H_{\infty, c}^p (X)$ est le groupe de cohomologie de De Rham des
$p$-formes diff{\'e}rentiables de classe $\mathcal{C}^{\infty}$ sur $X$ {\`a}
support compact et enfin $\Lambda^p (F)$ l'espace des $p$-formes
diff{\'e}rentielles de classe $\mathcal{C}^{\infty}$ sur $F$.

\begin{remark}
  
\end{remark}

On a aussi , $\bar{\partial}^2 = 0$ ; $\partial^2 = 0$.

Une forme diff{\'e}rentielle $\Omega \subset X$ d{\'e}finit un courant
not{\'e} $[f]$ sur $\Omega$ de la mani{\`e}re suivante :
\[ [f] (\varphi) = \int_{\Omega} f \wedge \varphi . \]
Il y a une application naturelle :
\[ H^{(p, q)} (\Omega) \tilde{\longrightarrow} H_{\tmop{cour}}^{(p, q)}
   (\Omega), \]
c'est l'isomorphisme de Dolbeault et une application naturelle :
\[ \check{H}^{(p, q)} (\Omega) \longrightarrow H_{\tmop{cour}}^{(p, q)}
   (\Omega) . \]

\section{R{\'e}solution de l'{\'e}quation $du = T$}

On va s'int{\'e}resser {\`a} une famille particuliere de domaines :

Soit $X$ une vari{\'e}t{\'e} diff{\'e}rentiable de dimension $n$ , $\Omega
\subset X$ un domaine born{\'e} {\`a} bord b$\Omega$ pour lequel on a $H^j  (b
\Omega) = 0$ pour $1 \leqslant j \leqslant n - 2$ o{\`u} $H^j  (b \Omega) =
j^{\tmop{iem}}$ groupe de cohomologie de De Rham .

On a $X = (X \setminus \Omega \nobracket) \cup \bar{\Omega}$ et $\Lambda^r
(F)$ est l'espace des $r$-formes diff{\'e}rentielles sur $F$.

On a la suite courte suivante :
\vspace{0,7em}

$0\longrightarrow \Lambda^{\bullet} (X) \longrightarrow \Lambda^{\bullet}  (X
\backslash \nobracket \Omega) \oplus \Lambda^{\bullet} ( \bar{\Omega})
\longrightarrow \Lambda^{\bullet}  (b \Omega) \longrightarrow 0$

\ \ \ \ \ \ $f$ \ \ $\longmapsto$ \ $(f_{|X \setminus \Omega}, f_{|
\bar{\Omega}})$ \ $\longmapsto$ \ $f_{|b \Omega}$

o{\`u} $\bullet = 0, 1, \ldots \nocomma, n$.

\vspace{0,7em}

On peut l'{\'e}crire en extension :

\vspace{0,7em}

0$\longrightarrow \Lambda^0 (X) \longrightarrow \Lambda^0  (X \backslash
\nobracket \Omega) \oplus \Lambda^0 ( \bar{\Omega}) \longrightarrow \Lambda^0 
(b \Omega) \longrightarrow 0$

\ \ \ \ \ $\downarrow d$ \ \ \ \ \ \ \ \ \ \ \ \ \ \ \ \ $\downarrow d$ \ \ \
\ \ \ $\downarrow d$ \ \ \ \ \ \ \ \ \ \ \ \ \ $\downarrow d$

0$\longrightarrow \Lambda^1 (X) \longrightarrow \Lambda^1  (X \backslash
\nobracket \Omega) \oplus \Lambda^1 ( \bar{\Omega}) \longrightarrow \Lambda^1 
(b \Omega) \longrightarrow 0$

. . . . . . . . . . . . . . . . . . . . . . . . . . . . . . .

.

. \ \ \ \ \ \ $\downarrow d$ \ \ \ \ \ \ \ \ \ \ \ \ \ \ \ $\downarrow d$ \ \
\ \ \ \ \ \ \ \ \ \ \ \ \ \ $\downarrow d$

0$\longrightarrow \Lambda^n (X) \longrightarrow \Lambda^n  (X \backslash
\nobracket \Omega) \oplus \Lambda^n ( \bar{\Omega}) \longrightarrow 0$

\vspace{0,7em}

On peut donc associer une suite longue sur le plan cohomologique suivante :

\vspace{0,7em}

$0 \longrightarrow H^0 (X) \longrightarrow H^0  (X \backslash \nobracket
\Omega) \oplus H^0 ( \bar{\Omega}) \longrightarrow H^0  (b \Omega)
\longrightarrow \ldots$

$\longrightarrow H^{n - 1} (X) \longrightarrow H^{n - 1}  (X \backslash
\nobracket \Omega) \oplus H^{n - 1} ( \bar{\Omega}) \longrightarrow H^{n - 1} 
(b \Omega) \longrightarrow$

$H^n (X) \longrightarrow H^n  (X \backslash \nobracket \Omega) \oplus H^n (
\bar{\Omega}) \longrightarrow 0$

\vspace{0,7em}

Tenant compte de $H^p ( b \Omega) = 0$ pour $p \neq 0$ et $p \neq n$ et $H^p (
\Omega) = 0$ pour $p \geqslant 1$,

on a $H^j  (b \Omega) = 0$ pour $j = 1, \ldots \nocomma \nocomma, n - 2$ \ par
cons{\'e}quent

$H^j (X) \backsimeq H^j  (X \setminus \Omega) \oplus H^j ( \bar{\Omega})$ pour
$1 < j \leqslant n - 2$ .

\begin{proposition}\label{P1}
  
\end{proposition}

Si $\Omega$ est contractile alors $H^j (X) = H^j  (X \setminus \Omega)$.

\begin{theorem}\label{th1}
  
\end{theorem}

Soit $X$ une vari{\'e}t{\'e} diff{\'e}rentiable de dimension $n$ pour laquelle
on a $H^j (X) = 0$ , $1 \leqslant j \leqslant n$. Alors pour tout domaine
contractile $\Omega$ et born{\'e} avec $H^j  (b \Omega) = 0$ , $1 \leqslant
j \leqslant n - 2$ , on a $\check{H}^j (\Omega) = 0$ , $1 \leqslant j
\leqslant n - 2$.

\vspace{0,7em}

Pour faire la preuve on a besoin du lemme suivant :

\begin{lemma}\label{lem1}
  
\end{lemma}

Sous les hypoth{\`e}ses de la proposition \ref{P1} , on a
\[ D^j ( \bar{\Omega}) \cap \ker d = d D^{j - 1} ( \bar{\Omega}) 
   \text{\tmop{pour}} 1 \leqslant j \leqslant n - 1. \]
\[  \]
\begin{proof}
  Soit $f \in D^{n - j} ( \bar{\Omega}) \cap \ker d$ $\Longrightarrow f \in
  D^{n - j} (X) \cap \ker d$ $\Longleftrightarrow [f] \in H^j (X)$, il existe
  $g$ une $(n - j - 1)$-forme diff{\'e}rentielle sur $X$ tel que $f = dg$ et
  puis que $f_{|X \setminus \Omega} = 0$ , alors $dg_{|X \setminus
  \bar{\Omega}} = 0$ d'o{\`u} $[dg_{|X \setminus \bar{\Omega}}] \in H^{n - j -
  1}  (X \setminus \bar{\Omega}) = H^{n - j - 1}  (X \setminus \Omega)$ donc
  il existe $u \in \Lambda^{n - j - 1}  (X \setminus \Omega)$ telle que $g_{|X
  \setminus \Omega} = u$.
  
  Soit $\tilde{u}$ une extension de $u$ {\`a} $X$. $v = g - d \tilde{u}$ est
  dans $\Lambda^{n - j - 1} (X) \Longrightarrow v = 0$ sur $X \setminus
  \Omega$ d'o{\`u} $v \in D^{n - j - 1} ( \bar{\Omega})$ et $dv = dg - d^2 
  \tilde{u} = f$ , donc $f = dv$ avec $v \in D^{n - j - 1} ( \bar{\Omega})$.
\end{proof}

On peut donc {\'e}tablir la preuve du th{\'e}or{\`e}me \ref{th1} :

\begin{proof}
  $\check{D}^j (\Omega) = (D^{n - j} ( \bar{\Omega}))'$ et $\check{D}^j
  (\Omega) \cap \ker d = d \check{D}^{j - 1} (\Omega)$ , $1 \leqslant j
  \leqslant n - 2$.
  
  Soit $C \in \check{D}^j (\Omega) \cap \ker d$ , $1 \leqslant j \leqslant n -
  2$.
  
  Posons $L_C : dD^{n - j} ( \bar{\Omega}) \longrightarrow \mathbb{C}$
  
  \ \ \ \ \ \ \ \ \ \ \ \ \ \ \ \ \ \ \ $d \varphi$ \ $\longmapsto$ \
  $\langle C \nocomma, \varphi \rangle$
  
  $L_C$ est bien d{\'e}finie car $D^{j - 1} ( \bar{\Omega}) = D^j (
  \bar{\Omega}) \cap \ker d$.
  
  Si $\varphi$ et $\varphi'$ sont deux $(n - p)$-forme diff{\'e}rentielle
  telles que $d \varphi = d \varphi'$ , alors
  \[ \varphi - \varphi' = d \theta \text{\tmop{et}}  \langle C, d \theta
     \rangle = \lim_{j \rightarrow + \infty} \langle C, d \theta_j \rangle = 0
     \Longrightarrow \langle C, \varphi \rangle = \langle C, \varphi' \rangle
     . \]
  $L_C$ est lin{\'e}aire car $d$ : $D^{n - j} ( \bar{\Omega}) \longrightarrow
  dD^{n - j} ( \bar{\Omega})$ est une application lineaire continue et
  surjective entre deux espaces de Fr{\'e}chet.
  
  Pour montrer que $L_C$ est continue , il suffit de voir que l'image
  r{\'e}ciproque d'un ouvert de $\mathbb{C}$ par $L_C$ est un ouvert . En
  effet on a $L_C od = C$ d'o{\`u} $L_C^{- 1} (U) = doC^{- 1} (U)$, par
  cons{\'e}quent on peut {\'e}tendre $L_C$ en un op{\'e}rateur lin{\'e}aire
  continu
  
  $\tilde{L}_C : D^{n - j + 1} ( \bar{\Omega}) \longrightarrow \mathbb{C}$ ,
  c'est un courant prolongeable et $d \tilde{L}_C \varphi = (- 1)^{n - j} C$
  et $\langle d \tilde{L}_C, \varphi \rangle = (- 1)^{n - j} \langle L_C, d
  \varphi \rangle = (- 1)^{n - j} \langle C, \varphi \rangle$. D'o{\`u} $S =
  (- 1)^{n - j} L_C$ est un courant prolongeable solution de $du = C$.
\end{proof}

\begin{remark}
  
\end{remark}

Le th{\'e}or{\`e}me reste vrai pour tout domaine contractile $\Omega$ {\`a}
bord v{\'e}rifiant $H^j  (b \Omega) = 0$ avec \ $1 \leqslant j \leqslant n -
2$. On ne sait pas cependant si pour un domaine contractile quelconque ce
th{\'e}or{\`e}me est vrai.

\section{R{\'e}solution du $\partial \bar{\partial}$ pour les courants
prolongeables}

Tenant compte du th{\'e}or{\`e}me \ref{th1} et des r{\'e}sultats de r{\'e}solution du
$\bar{\partial}$ pour les courants prolongeables obtenus par \cite{Samb} , on a le th{\'e}or{\`e}me suivant :

\begin{theorem}
  
\end{theorem}

Si $\Omega \Subset X$ est un domaine contractile strictement pseudoconvexe, alors pour tout

$(p, q)$-courant $d$-ferm{\'e} et prolongeable $\check{C}$ , il existe un $(p
- 1, q - 1)$-courant $\check{S}$ tel que $\partial \bar{\partial} \check{S} =
\check{C}$ sur $\Omega$ pour $2 \leqslant p + q \leqslant 2 n - 2$ .

Pour $\Omega \Subset X$ il existe une fonction $\varphi$ d{\'e}finie sur
$\bar{\Omega}$ strictement plurisousharmonique au voisinage $U_{\bar{\Omega}}$
de $\bar{\Omega}$ avec
\[ \Omega = \{z \in U_{\bar{\Omega}} : \varphi (z) < 0\} . \]
\begin{proof}
  Soit $C$ un $(p, q)$-courant , $1 \leqslant p \leqslant n - 1$ et $1
  \leqslant q \leqslant n - 1$ , $d$-ferm{\'e} d{\'e}fini sur $\Omega$ et
  prolongeable avec $2 \leqslant p + q \leqslant 2 n - 2$.
  
  Puisque d'apr{\`e}s le th{\'e}or{\`e}me \ref{th1} $\check{H}^{p + q} (\Omega) = 0$ ,
  il existe un courant prolongeable $h$ d{\'e}fini sur $\Omega$ tel que $dh =
  C$.
  
  $h$ est un $(p + q - 1)$-courant , il se d{\'e}compose en un $(p - 1,
  q)$-courant $h_1$ et en un $(p, q - 1)$-courant $h_2$. On a $dh = d (h_1 +
  h_2) = dh_1 + dh_2 = C$.
  
  Comme $d = \partial + \bar{\partial}$ , on a pour des raisons de bidegr{\'e}
  $\partial h_2 = 0$ et $\bar{\partial} h_1 = 0$ et d'apr{\`e}s le
  th{\'e}or{\`e}me principal dans \cite{Sam} , $h_1 =
  \bar{\partial} u_1$ et $h_2 = \partial u_2$ avec $u_1$ et $u_2$ des courants
  prolongeables d{\'e}finis sur $\Omega$.
  
  On a : $C = \partial h_1 \noplus + \bar{\partial} h_2$
  
  \ \ \ \ \ \ \ \ \ \ $= \partial \bar{\partial} u_1 + \bar{\partial}
  \partial u_2$
  
  \ \ \ \ \ \ \ \ \ \ $= \partial \bar{\partial}  (u_1 - u_2)$
  
  Posons $S = u_1 - u_2$ , $S$ est un $(p - 1, q - 1)$-courant prolongeable
  d{\'e}fini sur $\Omega$ tel que $\partial \bar{\partial} S = C$.
\end{proof}

\begin{remark}
  
\end{remark}

On peut continuer la d{\'e}composition du $(p + q - 1)$-courant $h$ en $( p -
2, q + 1)$ , $( p + 1, q - 2)$ ect. Mais pour des raison de bidegr{\'e}, on
constate que tous ces termes sont $d$-ferm{\'e}s. Par cons{\'e}quent, la
solution n'{\'e}tant pas unique, on a choisi celle qui convient.

\section{ Exemples de domaines $\Omega$ }

On prend $X =\mathbb{C}^n$. Alors l'{\'e}quation $( \ast)$ admet une solution
dans les cas suivants :

{\tmstrong{- La boule euclidienne. }}
\[ \Omega = B (o, r) = : \left\{ x = (x_1, \ldots, x_n) \in \mathbb{R}^n ;
   \sum_{j = 1}^n x_j^2 < 1 \right\} . \]
Ce cas est d{\'e}j{\`a} r{\'e}solut dans \cite{SBD}.

{\tmstrong{- Un pav{\`e}. }}
\[ \Omega =] a_1, b_1 [ \times \cdots \times] a_n, b_n [ \subset \mathbb{R}^n
   . \]
{\tmstrong{- Une hypersurface.}}
\[ \Omega =\mathbb{R}^n \times \{ y_{n + 1} > 0 \} \subset \mathbb{R}^{n +
   1} . \]
Cependant dans ce cas, $\Omega$ n'est pas born{\'e} et pour cela, on reviendra
plus tard sur ces domaines non born{\'e}s v{\'e}rifiant les hypoth{\`e}ses
cohomologiques du th{\'e}or{\`e}me.

\end{document}